\documentclass[12pt]{article}
\usepackage{amssymb,amsmath,epsfig,amsthm}

\newtheorem{theorem}{Theorem}[section]
\newtheorem{thm}[theorem]{Theorem}

\newtheorem{prop}[theorem]{Proposition}

\newtheorem{lemma}[theorem]{Lemma}
\newtheorem{cor}[theorem]{Corollary}
\newtheorem{corollary}[theorem]{Corollary}
\newtheorem{conjecture}[theorem]{Conjecture}
 
\theoremstyle{remark}
\newtheorem{remark}[theorem]{Remark}

\theoremstyle{definition}

\def\<{\langle}
\def\>{\rangle}
\def\Z{\mathbb{Z}}

\def\Q{\mathbb{Q}}

\def\G{\Gamma}

\def\o{\omega}
\def\a{\alpha}

\def\FP{{\rm{FP}}}

\title{Subgroups of direct products of two limit groups}
\author{Martin R. Bridson and James Howie}
\begin{document}
\maketitle

\begin{abstract} If  $\G_1$ and $\G_2$ are  limit
groups and $S\subset \G_1\times\G_2$   is of type $\FP_2$
 then $S$ has a subgroup of finite index that
is a product of at most two limit groups.
\end{abstract}

\section{Introduction}

In  \cite{BHMS} Miller, Short and the present authors
described the homological finiteness properties of subgroups
of direct products of free and surface groups. This completed a 
programme of research
that began with celebrated examples of Stallings \cite{stall} and Bieri \cite{bieri1}
and took an unexpected turn when Baumslag and Roseblade
\cite{BR} showed that the only finitely presented subdirect\footnote{Recall
that $S\subset A\times B$ is termed a subdirect product if its projection
to both $A$ and $B$ is surjective.}
products of two free groups are the obvious ones, namely free groups and
subgroups of finite index. 

This last 
theorem admits instructive proofs from several perspectives
\cite{BW}, \cite{short}, \cite{miller}; each of these proofs  lends
new insight and leads to an interesting generalisation, but each relies  heavily on the fact that one is dealing with
only two factors.  The search for a higher-dimensional
analogue culminated in the  Main Theorem of
\cite{BHMS}, which shows in particular that if $S$ is a subgroup
 of a direct product of $n$ free or surface groups and the
homology of $S$ (with trivial coefficients)
is finitely generated up to dimension $n$, 
then $S$ has a subgroup of finite index that is
isomorphic to a direct product of free or surface groups.

Upon the successful completion of his  project to characterise 
those finitely generated groups that have the same elementary or existential
theories as free groups (\cite{S1} 
to \cite{S6}; cf.~\cite{KM1, KM2}),
 Zlil Sela proposed a list of central questions concerning
these {\em limit groups} \cite{zlil:qu}. In particular he asked whether
the results of \cite{BHMS} could be extended to subdirect products
of limit groups. In the first two papers in this series \cite{bh1}, \cite{bh2},
we established results concerning the subgroup structure of limit groups
that enabled us to answer this question in the case where the limit groups have the
same elementary theory as a free group. In the present
article we turn our attention to the subdirect products of general limit groups.
We use our previous results in tandem with calcuations in group homology and
geometrical arguments based on Bass-Serre theory
 to prove the following analogue of the Baumslag-Roseblade theorem.

\begin{theorem}\label{main}
 Let $S\subset \G_1\times\G_2$ be a subdirect product
of  non-abelian limit groups.
 If $S$  is of type $\FP_2(\Q)$  
 and $L_i=S\cap \G_i$ is non-trivial for $i=1,2$, then
$S$ has finite index in $\G_1\times \G_2$. 
\end{theorem} 

Some simple reductions (\ref{reduce}) yield:

\begin{corollary}\label{fpres} Let $\G_1$ and $\G_2$ be  limit
groups.  If $S\subset \G_1\times\G_2$   is finitely presented
 then $S$ has a subgroup of finite index that
is a product of at most two limit groups.
\end{corollary}

This is a special case of the following result, which is proved in section \ref{cor}.

\begin{cor}\label{abfactors}
Let $\G_1,\dots,\G_n$ be limit groups, with $\G_3,\dots,\G_n$
abelian, and let $S\subset \G_1\times\cdots\times\G_n$ be a
subdirect product such that $S$ is of type $\FP_m(\Q)$,
where $m=\min(n,3)$.  Then $S$ is isomorphic to a subgroup
of finite index in a product of at most $m$ limit groups.
\end{cor}

For simplicity of exposition, the homology of a group $G$ in this
paper will always be with coefficients in a $\Q\,G$-module -- typically the
trivial module $\Q$.  
With minor modifications, our arguments also apply
with other coefficient modules, giving corresponding results under
the finiteness conditions $\FP_n(R)$ for other suitable rings $R$.

The final section of this paper contains some speculation concerning
the subgroups of the direct product of arbitrarily many limit groups.

\section{Limit groups, their decomposition and subgroup structure}

Since this is the third in a series of papers
(following \cite{bh1,bh2}), we shall only recall the minimal
necessary amount of information about limit groups. The reader unfamiliar
with this fascinating class of groups can consult the seminal papers of
Sela \cite{S1} to \cite{S6} or  the excellent introduction \cite{BF}.
 The papers of Kharlampovich and Myasnikov \cite{KM1,KM2} 
 approach the subject from a perspective more in keeping with
traditional combinatorial group theory, and a further
 perspective is developed in \cite{CG}.
 
\subsection{Limit groups}

A finitely generated group $\G$ is  a {\em limit group} if
for every finite set $T\subset\G$ there exists a homomorphism
$\phi$ from $\G$ to a free group so that $\phi|_T$ is injective.

Our results rely  on  the fact that limit groups
are the finitely generated subgroups of $\omega$-residually free tower ($\o$-rft)
 groups \cite{S2}.
Such tower groups are the fundamental groups of
tower spaces that are assembled from graphs, tori
and surfaces in a hierarchical manner. The number of
stages in the process of assembly is the {\em{height}}
of the tower. Each stage in the construction involves
the attachment of a surface along its boundary, or
the attachment of an $n$-torus $T$  along an
embedded circle representing a primitive element
of $\pi_1T$.

The {\em{height}} 
of a limit group $S$ is the minimal height of an 
$\o$-rft group
that has a subgroup isomorphic to $S$. Limit
groups of height $0$ are free products of 
finitely many free abelian groups (each of finite rank)
and surface
groups of Euler characteristic at most $-2$.

 The splitting described in the following proposition is obtained as follows:
embed $\G$ in an $\omega$-rft group $G$, take  the
graph of groups decomposition  that   the Seifert-van Kampen
Theorem associates to the addition of the final block in the tower, 
then apply Bass-Serre theory to get an induced graph of groups decomposition  of 
$\G$.

\begin{prop}\label{graph} If $\G$ is a limit group of
height $h\ge 1$, then $\G$ is the fundamental group of a
finite graph of groups with the following properties:
\begin{enumerate}
\item the vertex groups are  
limit groups of height $\le h-1$;
\item the edge groups are trivial or infinite-cyclic;
\item each edge is adjacent to a vertex where the vertex-group
is non-abelian.
\end{enumerate}
\end{prop}
 
Note also that any non-abelian limit group of height $0$  splits
as $A\ast_C B$ with $C$ infinite-cyclic or trivial.

\subsection{Homological finiteness implies finiteness}

The following is proved in \cite{bh1}.

\begin{thm}\label{1factor}
 Let $\G$ be a limit group and  $S\subset \G$  a subgroup.
 If $H_1(S,\Q)$ is finite dimensional
then  $S$ is finitely generated (and hence is a limit group).
\end{thm}

\subsection{Normal subgroups}

From \cite{bh1} we need:

\begin{thm}\label{norm} If $\G$ is a limit group
and $H\subset\G$ is a finitely generated subgroup,
then either $H$ has finite index in its normaliser
or else the normaliser of $H$ is abelian.
\end{thm}

\noindent{And from \cite{bh2} we need:}

\begin{theorem}\label{mh} Let $\G$ be a limit group that
splits as an amalgamated free product or HNN extension
with cyclic edge group $C\subset\G$.

If $N\subset\G$ is a non-trivial normal subgroup then there exists a
subgroup $\G_0\subset\G$ of finite index and an element $t\in N\cap \G_0$
such that  $\G_0$ is an HNN-extension with stable letter $t$ and
amalgamated subgroup $C\cap \G_0$.
\end{theorem}

\subsection{Double cosets}

We also require the following {\em Double Coset Lemma} from \cite{bh2}.

\begin{theorem}\label{dcl}
Suppose that, for $i=1,\dots,n$, $\G_i\ast_{A_i}$ is an HNN
extension with associated subgroup $A_i$ and stable letter $t_i$.
Let $\G=\G_1\times\cdots\times\G_n$, $A=A_1\times\cdots\times A_n$,
and let $S\subset\G$ be a subgroup containing $t_1,\dots,t_n$.
Then $H_n(S,\,\Q)$ contains a subgroup isomorphic to 
$$\Q\otimes_{\Q\,S}\Q\,\G\otimes_{\Q\,A}\Q,$$
 a rational vector space with basis the
set $S\backslash\G/A$ of double cosets $SgA$ with $g\in\G$.
\end{theorem}

In the proof of Theorem 1.1 we will need the above result in the
case $n=2$. In this case, if
$H_2(S,\,\Q)$ is finite dimensional then $|S\backslash\G/A|<\infty$.
Thus it is of interest to understand when a pair of subgroups
in a group  admit only finitely many double cosets. 
The following assertions are proved
in \cite{bh2}, and are  easily verified.

\begin{prop}\label{p:dc}
Let $G$ be a group, $A,B,C$ subgroups of $G$ with $A\subset C$
and $|B\backslash G/A|<\infty$, and $f:G\to H$ a homomorphism.
Then
\begin{enumerate}
\item $|(B\cap C)\backslash C/A|<\infty$, and
\item $|f(B)\backslash f(G)/f(A)|<\infty$.
\end{enumerate}
\end{prop}
 
\section{Some homological algebra}

In this section we gather the
homological algebra that we shall need in the proof of the main theorem.

\subsection{The LHS spectral sequence}\label{LHS}

Associated to any short exact sequence of groups $1\to A\to B\to C\to 1$
one has the Lyndon-Hochschild-Serre (LHS) spectral sequence in homology
(see \cite[p. 171]{brown},  or \cite[p.303]{hs}).
Given a $\Q B$-module $M$, the $E^2$ page of the spectral sequence has
$(p,q)$-entry $H_p(C, \, H_q(A,M))$, and the spectral sequence converges
to $H_*(B,M)$. 

The computations in this article concern $H_i(B)$ with trivial coefficients
and $i\le 2$. The terms of relevance on the $E^2$ page are
$$H_0(C, H_2(A)),\  H_1(C, H_1(A)),\ H_2(C, H_0(A)).$$

In each of our
calculations, the first and last of these terms will  be
finitely generated and we will be concerned with the issue 
of finite generation 
for $H_1(C, H_1(A))$. The only non-zero derivative hitting this group is
$d_2 : H_3(C, H_0(A))\to  H_1(C, H_1(A))$. Since the terms $E^\infty_{p,q}$
with $p+q=2$ are the composition factors of a filtration of $H_2(B)$, we
see that $H_2(B)$ will be finitely generated only if the
quotient of $H_1(C, H_1(A))$ by the image of $d_2$ is finitely generated.
But $H_3(C, H_0(A))$ will always be finitely generated in our setting,
and hence $H_2(B)$ will be finitely generated  only if $H_1(C, H_1(A))$
is finitely generated.

\subsection{Fox calculus}

Recall that a presentation $\langle X\mid R\rangle$ of a group $G$
gives rise to a partial free $\Q\,G$-resolution
$$F_2\overset{\partial_2}\to F_1\overset{\partial_1}\to\Q\,G\to\Q$$
of the trivial $G$-module $\Q$.  Here $F_1$ and $F_2$ are free modules
with bases in one-to-one correspondence with $X$ and $R$ respectively.

The boundary maps $\partial_1$ and $\partial_2$ are  easy
to calculate. $\partial_1$ is given by an $X\times 1$-matrix where
the $x$-entry is $1-x$, and $\partial_2$ is given by an
$R\times X$-matrix where the $(r,x)$-entry is the Fox derivative
$\partial r/\partial x$, defined inductively on the length of
$r$ by the rules:
$$\frac{\partial}{\partial x}(1) = 0 = \frac{\partial}{\partial x}(y) = \frac{\partial }{\partial x}(y^{-1}),~~y\in X\smallsetminus\{x\},$$
$$\frac{\partial }{\partial x}(x) = 1,~~\frac{\partial }{\partial x}(x^{-1}) = x^{-1},~~\frac{\partial }{\partial x}(uv) = \frac{\partial u}{\partial x} + u\frac{\partial v}{\partial x}.$$

\begin{lemma}\label{foxy}
Let $\G=G\ast_C$ be an HNN extension with stable letter $t$,
and let $M$ be a $\Q\,\G$-module on which $t$ acts trivially.
Then $H_1(\G,\,M)$ contains a subgroup that admits an
epimorphism onto $H_0(C,\,M)$.
\end{lemma}

\begin{proof}
Given a presentation $\langle X \mid  R\rangle$ of $G$, and a set
of generators $Y$ for $C$ (expressed as words in the alphabet $X$), 
we have a presentation for $\G$
of the form
$$\mathcal{P}=\langle X , t\mid  R, ytu_yt^{-1}~(y\in Y)\rangle$$
for $\G$ (where the $u_y$ are words in $X$).

We calculate $H_1(\G,\,M)$ by tensoring the resolution
associated to this presentation with $M$ and calculating $\partial_2$ using Fox
calculus. Specifically, we have
$$
\bigoplus_{y\in Y} M\oplus \bigoplus_{r\in R} M \overset{\partial_2}\to M\oplus \bigoplus_{x\in X} M  
\overset{\partial_1}\to M\to 0.
$$
With regard to the visible decompositions, $\partial_2$ is given by a matrix
of the form
$$\left(\begin{array}{cc}
\frac\partial{\partial t}(ytu_yt^{-1}) & \ast\\
0 & \ast\end{array}\right),$$
where $\frac\partial{\partial t}(ytu_yt^{-1}) = y-1$, 
while the restriction of $\partial_1$ to the first factor is multiplication
by $(1-t)$. But since $t$  acts trivially on $M$, 
multiplication by $(1-t)$ is the zero-map on $M$. Therefore the first factor $M$ in dimension 1 is contained in the group of $1$-cycles.
Let $A$ denote the subgroup of $H_1(\G,\, M)$ generated by this
group of $1$-cycles.

Now the first coordinate of the image of $\partial_2$ is contained
in the submodule $N$ of $M$ generated by $\{(y-1)M,~y\in Y\}$.  It follows
that the kernel of the natural epimorphism $M\to A$ is contained
in $N$. Thus there is an
epimorphism $A\to M/N$, and since $Y$ generates $C$, the 
quotient $M/N$ is the group of coinvariants $H_0(C,\, M)$.
\end{proof}

\subsection{A mapping cylinder argument}

We are grateful to Peter Kropholler for a helpful
conversation concerning the following lemma.

\begin{lemma}\label{MV}
If $G$ is of type $\FP_n(\Q)$ and acts cocompactly on 
a tree with all edge-stabilisers of type $\FP_n(\Q)$, then all vertex-stabilisers
are also of type $\FP_n(\Q)$.
\end{lemma}

\begin{proof}
If $V$, $E$ are the vertex and edge sets of the tree, then the 
hypothesis on edge-stabilisers
 is equivalent to $\Q E$ having a projective $\Q G$-resolution that is finitely generated
in dimensions $\le n$, and the assertion of the lemma
 is that the analogous property for $\Q V$ holds.

Since the action is cocompact, $\Q V$ is finitely generated as a $\Q
G$-module, so there is a projective $\Q G$-resolution $\mathcal F$
for $\Q V$ that is finitely generated in dimension $0$.

Suppose inductively that $k<n$ and there is a projective $\Q G$-resolution $\mathcal F$ for $\Q V$ that
is finitely generated
in dimensions $\le k$.
Let $\mathcal P$ be a projective $\Q G$-resolution for $\Q E$ that is 
finitely generated
in dimensions $\le n$.

Let $f:{\mathcal P}\to {\mathcal F}$ be a chain-map lifting the boundary map $\Q E\to \Q V$, and form
the mapping cylinder resolution ${\mathcal M}={\mathcal M}(f)$ of $\Q$.

  Thus ${\mathcal M}_i={\mathcal P}_{i-1}\oplus {\mathcal F}_i$
is finitely generated
in dimensions $i\le k$. The boundary maps $D_i:{\mathcal M}_i\to{\mathcal M}_{i-1}$ are defined by
$$D_i(x,y):=(d^{\mathcal P}(x),(-1)^{i-1}f(x)+d^{\mathcal F}(y)),$$
where $d^{\mathcal P}$ and $d^{\mathcal F}$ are the boundary maps in
${\mathcal P}$ and ${\mathcal F}$ respectively.

Since $k<n$ and $G$ is of type $\FP_n(\Q)$, the kernel of the
$k$-th boundary map $D_k:{\mathcal M}_k\to {\mathcal M}_{k-1}$ is finitely
generated.  Since also ${\mathcal P}_{k+1}$ is finitely
generated, there
is a finitely
generated $\Q G$-direct summand $L$ of $F_{k+1}$ such that
\begin{enumerate}
\item[(i)] $D_{k+1}({\mathcal P}_k\oplus L)=\mathrm{Ker}(D_k)=D_{k+1}({\mathcal M}_{k+1})$; and
\item[(ii)] $f_{k+1}({\mathcal P}_{k+1})\subset L$.
\end{enumerate}

Now suppose that $a\in \mathrm{Ker} (d_k:{\mathcal F}_k\to {\mathcal F}_{k-1})$.  Then $(0,a)\in \mathrm{Ker}(D_k)$,
so by (i) we can find $b\in {\mathcal P}_k$ and $c\in L$ such that 
$$(0,a) = D_{k+1}(b,c) = ( d_k(b) , (-1)^k f_k(b) + d_{k+1}(c) ).$$
Then $b \in \mathrm{Ker}(d_k) = \mathrm{Im}(d_{k+1})$, so we can find
$e\in {\mathcal P}_{k+1}$ such
that $b = d_{k+1}(e)$.  Then $f_k(b) = f_k(d_{k+1}(e)) = d_{k+1}(f_{k+1}(e))$,
and
$a = (-1)^k f_k(b) + d_{k+1}(c) = d_{k+1}( (-1)^k f_{k+1}(e) + c ) \in d_{k+1}(L)$ (using (ii)).

Hence $\mathrm{Ker}(d_k:F_k\to F_{k-1}) = d_{k+1}(L)$ is finitely generated.  It follows that there is a projective $\Q G$-resolution of
$\Q V$ that is finitely generated in dimensions $\le k+1$.

By induction, there is a projective $\Q G$-resolution of
$\Q V$ that is finitely generated in dimensions $\le n$.  In other words,
the vertex stabilisers of the action are all of type $\FP_n(\Q)$, as claimed

\end{proof}

\section{The main argument}

In this section we prove Theorem 1.1 and Corollary 1.2.

\subsection{Reduction to subdirect products} \label{reduce}

  Let $\G_1$ and $\G_2$ be  non-abelian limit groups. 
  Let $S\subset\G_1\times\G_2$ be
a subgroup and let $p_i:S\to\G_i$
be the natural projection. If $S$ is of type $\FP_2(\Q)$, then
 $H_1(S,\,\Q)$ is finite dimensional, so
$H_1(p_i(S),\,\Q)$ is finite dimensional for $i=1,2$.
By Theorem \ref{1factor}, $p_i(S)\subset\G_i$ is a limit group.
Replacing $\G_i$ by $p_i(S)$, we reduce to the case where
 $S$ projects onto each factor $\G_i$.

Let $L_i=S\cap \G_i$. If $L_1$ were trivial, then $S\cong p_2(S)$
would be a limit group, satisfying the conclusion of Corollary \ref{fpres}.
Thus we may assume that $L_i$ is non-trivial for $i=1,2$
(and Corollary \ref{fpres} becomes a weak form of Theorem \ref{main}).

\begin{remark}\label{normal}
Since $p_i(S)=\G_i$ and  $L_i$ is normal in $S$,
we also have that $L_i=p_i(L_i)$
is normal in $\G_i$.
\end{remark}

\subsection{Subgroups of finite index} 

Throughout the proof we shall repeatedly pass to subgroups of finite
index in one or other of the $\G_i$. When we do so, we shall be assuming
that $S$ is replaced with the inverse image of this subgroup under the
projection $p_i:S\to\G_i$ and $\G_j$ ($j\neq i$) is replaced
by $p_jp_i^{-1}(H_i)$. This does not affect $L_j$.

\begin{remark}\label{fi}
 It is because of the need to pass to subgroups of 
finite index that we have assumed that $S$ is of type
$\FP_2$ rather than simply that $H_2(S)$ is finitely generated:
the former hypothesis passes to subgroups of finite index, whereas
the latter does not necessarily. 
\end{remark}

\subsection{Useful lemmas and special cases}

Much
of our proof will revolve around the LHS spectral sequences associated
to the short exact sequences for the projections of $S$ to the two factors  
 $1\to L_2\to S\to\G_1\to 1$ and $1\to L_1\to S\to\G_2\to 1$. In particular,
following the discussion in subsection \ref{LHS}, we shall use the hypothesis
that $S$ is of type $\FP_2$ by assuming that
$H_1(\G_i,\, H_1(L_j,\,\Q))$ is finitely generated for $j\neq i$. In the ensuing argument
we shall make several appeals to the following:

\begin{lemma}\label{ll} Let $S$ be a subdirect product of two
non-abelian limit groups $\G_1$ and $\G_2$, such that the intersections
$L_1=S\cap\G_1$ and $L_2=S\cap\G_2$ are both non-trivial.
If either 
$H_1(L_1,\,\Q)$ or $H_1(L_2,\,\Q)$
is finite dimensional, then $S$ has finite index in $\G_1\times\G_2$.
\end{lemma}

\begin{proof} Suppose that $H_1(L_1,\,\Q)$ is finite dimensional.
It follows from Theorem \ref{1factor} that $L_1$
is finitely generated, and from Theorem \ref{norm} and Remark
\ref{normal} that $L_1$ has finite
index in $\G_1$. Therefore, after replacing $\G_1$ by $L_1$,
we can split $1\to L_2\to S\to \G_1\to 1$  to obtain
$S= \G_1\times L_2$. Thus
$L_2$, being a retract of $S$, is finitely generated, so of finite index
in $\G_2$.
\end{proof}

\begin{prop}\label{prop}
Let $G_1$ be a group and let $\G_2$ be a limit group. 
Let $S\subset G_1\times\G_2$
be a subgroup of type $\FP_2(\Q)$ that projects onto $\G_2$, and
that intersects $\G_2$ non-trivially.
Suppose that $\G_2$ can be expressed as an amalgamated free
product or HNN extension with cyclic edge-group $C$.

If $C\cap S$ has finite index in
$C$, then  $H_1(G_1\cap S,\,\Q)$ is finite dimensional.
\end{prop}

\begin{proof} Let $L_1=G_1\cap S$ and let $L_2=\G_2\cap S$.
 Since the cyclic subgroups of 
 $\G_2$ are closed in the profinite topology,
  we may pass to a subgroup of finite
index which intersects $C$ in $C\cap L_2$ and hence assume
that $C\subset L_2$.

By appealing to Theorem \ref{mh} and passing to a further
subgroup of finite index,
we may assume that $\G_2$ is an HNN extension $G\ast_C$ with stable
letter $t\in L_2$. 

Consider the Lyndon-Hochschild-Serre
spectral sequence in $\Q$-homology corresponding to the group
extension
$$1\to L_1\to S\to\G_2\to 1.$$

The action of $S$ by conjugation on $L_1$ induces
an action of $\G_2=S/L_1$ on $H_*(L_1,\,\Q)$ in which
$L_2$ acts trivially. In particular $t$ and $C$ act trivially on
the $\G_2$-module $H_k(L_1,\,\Q)$, for each $k$.

The $(1,1)$-term $E_{1,1}^3$ on the $E^3$-page of this spectral
sequence survives to infinity, and so is isomorphic to a
subquotient of $H_2(S,\,\Q)$, which is finite-dimensional, by
hypothesis.  Hence $E_{1,1}^3$ is finite-dimensional.  But
$E_{1,1}^3$ is the cokernel of the differential $d:E_{3,0}^2\to E_{1,1}^2$
from the $E^2$-page. Moreover, $E_{3,0}^2=H_3(\G_2,\, H_0(L_1,\,\Q))=H_3(\G_2,\,\Q)$ is finite-dimensional, since the limit
group $\G_2$ has type $\FP_\infty(\Q)$.  Hence 
$E_{1,1}^2=H_1(\G_2,\, H_1(L_1,\,\Q))$ is finite-dimensional.

Since $t$ acts trivially on $H_*(L_1,\,\Q)$, we may apply
Lemma \ref{foxy}, with $\G=\G_2$ and $M=H_1(L_1,\,\Q)$.
It follows that $H_0(C,\, H_1(L_1,\,\Q))$
is finite-dimensional. And $H_0(C,\, H_1(L_1,\,\Q))\cong H_1(L_1,\,\Q)$
 since $C$ acts trivially on $H_*(L_1,\,\Q)$.
\end{proof}

The spectral sequence argument ending with the penultimate paragraph
of the preceding proof admits a straightforward generalisation as 
follows.  If $S$ is of type $\FP_n$, $k\le n-1$ and
 $H_j(L_1,\,\Q)$ is finite dimensional  for $j < k$, then
$E_{1,k}^2=H_1(\G_2,\, H_k(L_1,\,\Q))$ is finite-dimensional. Taking
$M= H_k(L_1,\,\Q)$ in the final paragraph then implies
$H_k(L_1,\,\Q)\cong H_0(C,\,H_k(L_1,\,\Q))$ is finite-dimensional.

To generalise further, recall that
Theorem 2.4 is a weak form of what we actually proved in \cite{bh2}
concerning the existence of stable letters in normal subgroups. 
Correspondingly, we can weaken the hypotheses on $\G_2$ and $C$ in
 the above proposition, and appeal to 
  \cite[Corollaries 2.2 and 3.2]{bh2} instead of Theorem 2.4. These
  remarks lead to the following generalisation (which will not be needed
  in the sequel).

 \begin{prop} \label{genProp} Let $G_1$ and $G_2$ be groups.
Suppose $G_2$ is of type $\FP_{n+1}(\Q)$ and admits an 
$l$-acylindrical\footnote{A splitting is $l$-acylindrical
 if arcs of length $l$ in the 
 associated Bass-Serre tree have trivial stabilizers.} splitting
 (for some $l>0$) over a subgroup
$C$ that is closed in the profinite topology.

Let $S\subset G_1\times G_2$ be a subgroup of type $\FP_{n}$ that 
intersects $G_2$ non-trivially and intersects $C$ in a subgroup of
finite index. Then $H_i(G_1\cap S, \, \Q)$ is finite dimensional for
$i\le n-1$.
\end{prop}

In the light of Lemma 4.3, the following is an immediate
consequence of Proposition \ref{prop}.

\begin{cor}\label{entriv}
Let $\G_1$ and $\G_2$ be non-abelian limit groups.
Let $S\subset \G_1\times\G_2$
be a subdirect product of type $\FP_2(\Q)$ 
 that intersects each of the factors non-trivially.
Suppose that $\G_2$ has a graph-of-groups decomposition
with a cyclic edge-group $C$.

If $C\cap S$ has finite index in
$C$, then $S$ has finite index in $\G_1\times\G_2$.
\end{cor}

In the case where $C$ is trivial we have:

\begin{cor}\label{freeprod} Let  $\G_1$ and $\G_2$ be non-abelian
limit groups and let  $S\subset\G_1\times\G_2$ be a subdirect product
of type $\FP_2(\Q)$ that intersects each of the factors non-trivially.
 If $\G_2$ is a nontrivial free product, then
$S$ has finite index in $\G_1\times\G_2$.
\end{cor}

The remaining non-abelian limit groups of height $0$ are
covered by:

\begin{cor}\label{surface}  Let  $\G_1$ and $\G_2$ be non-abelian 
limit groups with $\G_2$ the fundamental
group of a compact surface. Let $S\subset\G_1\times\G_2$ be a 
subdirect product  that intersects each of the factors non-trivially.

If $S$ is of type $\FP_2(\Q)$, then
$S$ has finite index in $\G_1\times\G_2$.
\end{cor}

\begin{proof}
Passing to an index $2$ subgroup of $\G_2$ if necessary, we
may assume that the surface in question is orientable, and by
Corollary \ref{freeprod} we may assume that it is  closed.

Fix $\a\ne 1$ in $L_2$.
By a theorem of Scott \cite{scott}, there is a finite-sheeted cover
of $\Sigma$ in which $\a$ lifts to a simple closed curve.
Replacing $\G_2$ by the corresponding subgroup of finite index, we
may asume that $\G_2$ splits over the cyclic subgroup $\langle\a\rangle\subset L_2$.  The result then follows from Corollary \ref{entriv}.
\end{proof}

\medskip

The following observation will be used to control the way  in which
$S$ intersects each non-abelian vertex group
in the decomposition of $\G_2$ provided by Proposition \ref{graph}.

\begin{lemma}\label{l:vtriv}
Let $\G_1$ and $\G_2$ be non-abelian
limit groups and $S\subset\G_1\times\G_2$
a subdirect product of type $\FP_2(\Q)$ with $L_i=\G_i\cap S$
non-trivial for $i=1,2$. Consider a
graph-of-groups decomposition $\mathcal G$ of $\G_2$ that has
cyclic edge groups and a vertex group $\G_v$ that
is a non-abelian limit group.
Then $\G_v\cap L_2\ne\{1\}$.
\end{lemma}

\begin{proof} In $\mathcal G$
we fix an edge $e$  incident to the vertex $v$ and denote the corresponding
edge group   $\G_e\subset\G_2$.
We also fix a graph-of-groups decomposition of $\G_1$ that
has a cyclic edge group $\G_f$ (see Proposition \ref{graph} and the remark following it).

Using Theorem \ref{mh}, we may replace $\G_1$
and $\G_2$ by subgroups of finite index 
and assume that they split as HNN extensions
 with amalgamated subgroups $\G_e, \G_f$ and
 stable letters in $L_1$ and $L_2$, respectively. 
The Double Coset Lemma (Theorem \ref{dcl})
then tells us that $H_2(S,\,\Q)$ contains a subgroup
isomorphic to 
$$\Q(S\backslash (\G_1\times\G_2)/(\G_f\times\G_e)).$$ 

Since $H_2(S,\,\Q)$ is finitely generated, 
$|S\backslash (\G_1\times\G_2)/(\G_f\times\G_e)|<\infty$. So,
by Proposition \ref{p:dc}, $|(S\cap (\G_f\times\G_2))\backslash (\G_f\times\G_2)/(\G_f\times\G_e)|<\infty$.
Applying the projection $p_2:S\to\G_2$, we have 
$|p_2(S\cap (\G_f\times\G_2))\backslash\G_2/\G_e|<\infty$, by Proposition \ref{p:dc}.
Hence
$$|\G_v\cap p_2(S\cap (\G_f\times\G_2))\backslash\G_v/\G_e|<\infty,$$
again by Proposition \ref{p:dc}.
Now, $S\cap (\G_f\times\G_2)$ is generated by $L_2=p_2(L_2)$
and a cyclic group, so if $\G_v\cap L_2=\{1\}$ then
$C=\G_v\cap p_2(S\cap (\G_f\times\G_2))$ would
be a cyclic group.

But if $C$ were cyclic, then in the non-abelian, residually-free group
$\G_v$ we would have two cyclic subgroups, $C$ and $\G_e$, such that $|C\backslash\G_v/\G_e|<\infty$. But this is impossible,
for  there is   an epimorphism, $\phi$ say, from
$\G_v$ onto a non-abelian free group $F$, where the conclusion
(from Proposition \ref{p:dc}) that the cyclic subgroups $\phi(C)$ and $\phi(\G_e)$
satisfy $|\phi(C)\backslash F/\phi(\G_e)|<\infty$ is absurd.
\end{proof}

\subsection{Completing the proof of Theorem \ref{main}}

We have non-abelian limit groups $\G_1,\G_2$ and a subdirect
product $S\subset \G_1\times\G_2$ of type $\FP_2(\Q)$
such that $L_1=S\cap\G_1$
and $L_2=S\cap\G_2$ are non-trivial.   
We want to prove that $S$ has finite index in
$\G_1\times\G_2$. We shall do so by induction on
$$h:=\rm{height}(\G_1) + \rm{height}(\G_2),$$
 where {\em height}
is defined as in subsection 2.1. Let $h_i=\rm{height}(\G_i)$ for
$i=1,2$.

Corollaries \ref{freeprod} and \ref{surface} cover the case where
$h_2=0$ and the case where $\G_2$ is freely decomposable.
Thus, by Proposition \ref{graph}, we may assume that
$\G_2$ has a graph-of-groups decomposition 
$\G_2\cong\pi_1({\mathcal G},X)$ 
with infinite-cyclic
edge groups such that each edge is incident at a vertex where
the vertex group  is a non-abelian limit group of height $\le h_2-1$.
 Let $\G_v$ be such a vertex group and let   $\G_e$ be an
incident edge group.
 Lemma \ref{l:vtriv} tells us that $\G_v\cap L_2\ne\{1\}$.
Corollary \ref{entriv}  tells us that $S$ has finite index in $\G_1\times\G_2$
unless $\G_e\cap L_2=\{1\}$. Thus we proceed with the assumption
that $\G_e\cap L_2=\{1\}$ for all edges $e$ in the graph $X$. (This
will ultimately lead to a contradiction.)

Let $p_2:S\to\G_2$ be the natural projection and let $S_e=p_2^{-1}(\G_e)$.

\smallskip
\noindent{\em Claim:} $S_e$ is a limit group.

\smallskip
 Since $\G_e\cap L_2$ is trivial, 
the projection   $S\to \G_1$, which has kernel $L_2$, restricts to an injection
on $S_e$. Thus $S_e$ is isomorphic to a subgroup of a limit group
and the claim will follow from Theorem \ref{1factor} if we  prove 
that $H_1(S_e,\,\Q)$ is finite dimensional.

By Theorem \ref{mh}, we may
assume that 
 $\G_2$ splits as an HNN extension
with amalgamated subgroup  $\G_e$ 
and stable letter $t\in L_2$. (This reduction involves passing to a
subgroup of finite index  $S_e'\subset S_e$, but if
$H_1(S_e',\,\Q)$ is finite dimensional, then  $H_1(S_e,\,\Q)$ is.)

Let $M=H_1(L_1,\,\Q)$. Analysing the LHS spectral sequence for
$$1\to L_1\to S\overset{p_2}\to \G_2\to 1$$ 
 exactly
as in the proof of Proposition \ref{prop},
we see that since $H_2(S,\,\Q)$ is finite dimensional,  
$H_1(\G_2,\, M)$ is finite dimensional, and hence so is
$H_0(\G_e,\,M)$, by Lemma \ref{foxy}.

We can now calculate $H_1(S_e,\,\Q)$ using the LHS 
spectral sequence for 
$$1\to L_1\to S_e\overset{p_2}\to \G_e\to 1.$$ 
The relevant terms on the $E^2$ page are $H_1(\G_e,H_0(L_1,\,\Q))=H_1(\G_e,\,\Q)\cong\Q$
and $H_0(\G_e, \,M)$, which we have just argued is finitely generated. Thus
$H_1(S_e,\,\Q)$ is finite dimensional and our claim is proved.

\medskip
 
The graph of groups decomposition $\G_2\cong\pi_1({\mathcal G},X)$
pulls back to a graph of groups decomposition 
$S\cong\pi_1({\mathcal S},X)$   with the same underlying
graph $X$. Here the vertex and edge
groups have the form $S_v=p_2^{-1}(\G_v)$ and $S_e=p_2^{-1}(\G_e)$ respectively,
as $v$ and $e$ vary over the vertices and edges of $X$. 
We have just shown that the edge groups $S_e$ are all
limit
groups. In particular each $S_e$ has a finite classifying space
\cite{AB}, and so is of type $\FP_\infty(\Q)$.

We are assuming that $S$ is of type $\FP_2(\Q)$, so by
Lemma \ref{MV}, we can now deduce that {\em the vertex groups
$S_v$ are all of type $\FP_2(\Q)$. }

\smallskip
For the remainder of the proof we
focus our attention on a single non-abelian vertex
group $V:=\G_v$ in $\mathcal G$ and the corresponding vertex group
 $S_v$ in $\mathcal S$. 
Let $\Lambda$ be the
projection of $S_v$ to $\G_1$. Since $H_1(S_v,\,\Q)$ is finite
dimensional, so is $H_1(\Lambda,\,\Q)$. By Theorem \ref{1factor},
$\Lambda$ is a limit
group (of  height at most $h_1$). Moreover, $V$ is a limit
group of height strictly less than $h_2$ (see \ref{graph}).

Thus we have a product of limit groups $\Lambda\times V$
with the sum of heights $<h$, and a subgroup $S_v\subset \Lambda\times V$
that is
of type $\FP_2(\Q)$ and projects onto each factor. 
Moreover $S_v\cap\Lambda =L_1$, which is non-abelian since
it is non-trivial and normal in the non-abelian limit group $\G_1$. 
Hence, in particular, $\Lambda$ is non-abelian.

We also have, by hypothesis, that $V$ is non-abelian and that
$V\cap L_2\ne\{1\}$. Thus we may apply our inductive hypothesis
to conclude that $S_v$ has finite index in $\Lambda\times V$. Since $L_1=S_v\cap \Lambda$,
we conclude that $L_1$ has finite index in $\Lambda$ and hence is
finitely generated. Lemma \ref{ll} then
implies that $S$ has finite index in $\G_1\times\G_2$. This
contradicts the fact that $S\cap \G_e$ was assumed to
be trivial for each of the infinite-cyclic edge groups $\G_e\subset
\G_2$. \hfill $\square$

\section{Abelian factors}\label{cor}

Our main theorem concerns a subdirect product of two non-abelian
limit groups.  The result is false as stated in the case where
one or both of the limit groups is allowed to be abelian.

An abelian limit group is just a free abelian group of finite rank.
In particular, {\em every} subgroup of a direct product
of two such groups is of type $\FP_\infty(\Q)$.  More generally,
one can easily construct counterexamples to the exact analologue
of the main theorem in which one of the factors is allowed to
be non-abelian, along the following lines.

Let $\G_1$ be a non-abelian limit group, and $f_1:\G_1\to\Z^k$ an epimorphism for some $k\ge 1$.  Let $\G_2=\Z^m$ for some
$m>k$, and let $f_2:\G_2\to\Z^k$ be an epimorphism.  Let 
$S\subset \G_1\times\G_2$ be the
fibre product of $f_1$ and $f_2$. Then $S\cong\G_1\times\Z^{m-k}$
is a subdirect product of type $\FP_\infty$ and $S\cap \G_i$
is non-trivial for $i=1,2$, but $S$ does not have finite index
in $\G_1\times\G_2$.

It is not difficult to show that {\em all} counterexamples to the
exact analogue of the main theorem are, up to finite index,
isomorphic to one of the above
types.  In particular, in each case the group $S$ under consideration
is itself a direct product of limit groups, up to finite index.  

This observation is a special case of Corollary \ref{abfactors}.
We are now ready to prove Corollary \ref{abfactors}, which we restate
here as

\begin{cor}
Let $\G_1,\dots,\G_n$ be limit groups, with $\G_3,\dots,\G_n$
abelian, and let $S\subset \G_1\times\cdots\times\G_n$ be a
subdirect product such that $S$ is of type $\FP_m(\Q)$
(where $m=\min(n,3)$).  Then $S$ is isomorphic to a subgroup
of finite index in a product of at most $m$ limit groups.
\end{cor}

\begin{proof}
If all the $\G_i$ are abelian, then $S$ is a free abelian
group of finite rank, and hence a limit group. 
If, on the other hand, none of the $\G_i$ are abelian, then the result
follows from Corollary \ref{fpres}. Hence we may
assume that $n\ge 2$, that $\G_1$ is non-abelian, and that $\G_n$
is abelian.

If $\G_{n-1}$ is abelian, let $G_{n-1}$ be the image of $S$
under the direct projection
$\G_1\times\cdots\times\G_n\to\G_{n-1}\times\G_n$.  The result
for $S\subset \G_1,\dots,\G_n$
follows easily from the corresonding result for
$S\subset\G_1\times\cdots\times\G_{n-2}\times G_{n-1}$.

Hence we may assume that there is precisely one abelian factor,
namely $\G_n$,
and hence that $n\in\{2,3\}$.  Let $L_n=S\cap\G_n$.
Since $\G_n$ is abelian, we may assume (up to replacing
$\G_n$ by a finite index subgroup) that $L_n$ is a direct
factor of $\G_n$.  Say $\G_n\cong L_n\oplus M$.

Let $\pi:\G_n\to L_n$ be the projection with kernel $M$, and
define $\theta:=\pi\circ p_n:\G_1\times\cdots\times\G_n\to L_n$.
Then $\theta|_S$ is a splitting of the inclusion map
$L_n\to S$, so $S\cong T\times L_n$, where
$T=\mathrm{Ker}(\theta|_S)=S\cap(\G_1\times\cdots\times\G_{n-1}\times M)$.

Now $T\cap M=\{1\}$, and so $T$ is isomorphic to a subgroup of
$\G_1\times\cdots\times\G_{n-1}$.  Moreover, $T$ is of type $\FP_m(\Q)$
since $S\cong T\times L_n$ is of type $\FP_m(\Q)$.  It follows
from the main theorem
that $T$ has finite index in $\G_1\times\cdots\times\G_{n-1}$.

Hence $S$ is, up to finite index, a product of at most $m$ limit
groups, as claimed.
\end{proof}

\section{The natural conjecture for $n$ factors}

The results of this paper, and of \cite{bh2}, encourage us to make
the following conjecture, which would give a positive answer to a
question of Sela \cite{zlil:qu}.

\begin{conjecture}
Let $\G_1,\dots,\G_n$ be limit groups, and $S$ a subgroup of
$\G_1\times\cdots\times\G_n$ such that $S$ has type $FP_n(\Q)$.
Then $S$ has a subgroup of finite index that is itself a direct
product of $n$ or fewer limit groups.
\end{conjecture}

More precisely, we believe:

\begin{conjecture}
Let $\G_1,\dots,\G_n$ be non-abelian limit groups and  $S\subseteq
\G_1\times\dots\times\G_n$  a subdirect product 
that intersects each $\G_i$ non-trivially. If
exactly $k>0$ of the intersections $S\cap\G_i$ are not finitely
generated, then
 there exists a subgroup of finite index $S_0\subset S$
such that $H_k(S_0,\,\Q)$ is infinite dimensional.
\end{conjecture}

The arguments in \cite{bh1} and
\cite{bh2}  prove these conjectures in the case where the $\G_i$
are subgroups of elementarily free groups.  Those arguments used 
in a crucial way the
fact that the top of the $\o$-residually free tower
decomposition is a quadratic block.  It does not seem plausible to
extend the methods of \cite{bh2} to the case where each of the
$\G_i$ has an abelian block at the top of its $\o$-rft decomposition.

In the present paper, we addressed the case where at most $2$ of the
$\G_i$ are non-abelian, which reduces easily to the case  
of precisely two factors, each non-abelian. This more algebraic
approach (which is closer to the argument we used with
Miller and Short in the surface case \cite{BHMS}) may eventually
lead to a proof in the general case. However, there are considerable
difficulties associated to the fact that one must repeatedly pass
to subgroups of finite index. It seems possible that these
may be overcome by converting
many of our arguments based on the finiteness of various homology  groups
into arguments based on properties of corresponding homology
functors, but we do not know how to do this at present.

\medskip\centerline{\bf Authors' addresses}

\smallskip\begin{center}\begin{tabular}{lll}
Martin R. Bridson &\qquad\qquad & James Howie\\
Department of Mathematics && Department of Mathematics\\
Imperial College London && Heriot-Watt University\\
London SW7 2AZ && Edinburgh EH14 4AS\\
{\tt m.bridson@imperial.ac.uk} && {\tt J.Howie@hw.ac.uk}
\end{tabular}\end{center}

\end{document}